\documentclass[12pt]{amsart}
\usepackage{amsmath,amssymb,amsthm,amsfonts}
\usepackage{setspace,xypic}
\usepackage{graphicx}

\newlength{\figboxwidth}
\newcommand\inv{^{-1}}

\newcommand\G{\Gamma}

\newcommand\g{\gamma}
\newcommand\Ra{\mathbb R}

\newcommand\Za{\mathbb Z}

\newcommand\Ha{\mathbb H}
\mathchardef\GG="321D

 \DeclareMathOperator{\Isom}{Isom}

\DeclareMathOperator{\Sol}{Sol}
\DeclareMathOperator{\SOL}{Sol}
\DeclareMathOperator{\QI}{QI}\DeclareMathOperator{\Bilip}{Bilip}

\newtheorem{theorem}{Theorem}[section]

\newtheorem{lemma}[theorem]{Lemma}

\newtheorem{defn}[theorem]{Definition}

\newtheorem{question}[theorem]{Question}
\newtheorem{conjecture}[theorem]{Conjecture}

\newtheorem{step}{Step}

\begin{document}

\title[Quasi-isometries of solvable groups]{Quasi-isometries and rigidity of
solvable groups}
\author{Alex Eskin, David Fisher and Kevin Whyte}
\thanks{First author partially supported by NSF grant DMS-0244542.
Second author partially supported by NSF grants DMS-0226121 and
DMS-0541917. Third author partially supported by NSF grant
DMS-0349290 and a Sloan Foundation Fellowship.}

\begin{abstract}
In this note, we announce the first results on quasi-isometric
rigidity of non-nilpotent polycyclic groups.  In particular, we
prove that any group quasi-isometric to the three dimenionsional
solvable Lie group $\SOL$ is virtually a lattice in $\SOL$.  We
prove analogous results for groups quasi-isometric to
$\Ra{\ltimes}\Ra^n$ where the semidirect product is defined by a
diagonalizable matrix of determinant one with no eigenvalues on
the unit circle. Our approach to these problems is to first
classify all self quasi-isometries of the solvable Lie group.  Our
classification of self quasi-isometries for $\Ra{\ltimes}\Ra^n$
proves a conjecture made by Farb and Mosher in \cite{FM3}.

Our techniques for studying quasi-isometries extend to some other
classes of groups and spaces. In particular, we characterize groups
quasi-isometric to any lamplighter group, answering a question of de
la Harpe \cite{dlH}. Also, we prove that certain Diestel-Leader
graphs are not quasi-isometric to any finitely generated group,
verifying a conjecture of Diestel and Leader from \cite{DL} and
answering a question of Woess from \cite{SW,Wo}.  We also prove that
certain non-unimodular, non-hyperbolic solvable Lie groups are not
quasi-isometric to finitely generated groups.

The results in this paper are contributions to Gromov's program
for classifying finitely generated groups up to quasi-isometry
\cite{Gr2}.  We introduce a new technique for studying
quasi-isometries, which we refer to as {\em coarse
differentiation}.\\

{\it Dedicated to Gregory Margulis on the occasion of his $60$th
birthday.}
\end{abstract}

\maketitle

\section{Introduction and statements of rigidity results}
\label{section:rigidity}

For any group $\G$ generated by a subset $S$ one has the
associated Cayley graph, $C_{\G}(S)$.  This is the graph with
vertex set $\G$ and edges connecting any pair of elements which
differ by right multiplication by a generator.   There is a
natural $\G$ action on $C_{\G}(S)$ by left translation.   By
giving every edge length one, the Cayley graph can be made into a
(geodesic) metric space.   The distance on $\G$ viewed as
the vertices of the Cayley graph is the {\em word metric}, defined
via the norm:

 $$\|\g\|=\inf\{\text{length of a word in the generators }S
    \text{ representing } \g \text{ in } \G.\}$$

Different sets of generators give rise to different metrics and
Cayley graphs for a group but one wants these to be equivalent.  The
natural notion of equivalence in this category is {\em quasi-isometry}:

\begin{defn}
\label{defn:qi} Let $(X,d_X)$ and $(Y,d_Y)$ be metric spaces. Given
real numbers $K{\geq}1$ and $C{\geq}0$,a map $f:X{\rightarrow}Y$ is
called a {\em $(K,C)$-quasi-isometry} if
\begin{enumerate} \item
$\frac{1}{K}d_X(x_1,x_2)-C{\leq}d_Y(f(x_1),f(x_2)){\leq}Kd_X(x_1,x_2)+C$
for all $x_1$ and $x_2$ in $X$, and, \item the $C$ neighborhood of
$f(X)$ is all of $Y$.
\end{enumerate}
\end{defn}

If $\G$ is a finitely generated group, $\G$ is canonically
quasi-isometric to any finite index subgroup $\G'$ in $\G$ and to
any quotient $\G'' = \G / F$ for any finite normal subgroup $F$.
The equivalence relation generated by these (trivial)
quasi-isometries is called {\em weak commensurability}.   A group
is said to {\em virtually} have a property if some weakly
commensurable group does.

In his ICM address in 1983, Gromov proposed a broad program for
studying finitely generated groups as geometric objects, \cite{Gr2}.
Though there are many aspects to this program (see \cite{Gr3} for a
discussion), the principal question is the classification of finitely
generated groups up to quasi-isometry.  By construction, any
finitely generated group $\G$ is quasi-isometric to any space on which
$\G$ acts properly discontinuously and cocompactly by isometries. For
example, the fundamental group of a compact manifold is
quasi-isometric to the universal cover of the manifold (this is called
the Milnor-Svarc lemma).  In
particular, any two cocompact lattices in the same Lie group $G$ are
quasi-isometric.  One important aspect of Gromov's
program is that it allows one to generalize many invariants,
techniques, and questions from the study of lattices to all finitely
generated groups.

Given the motivations coming from the study of lattices, one of
the first questions in the field is whether a group
quasi-isometric to a lattice is itself a lattice, at least
virtually.  This question has been studied extensively.  For
lattices in semisimple groups this has been proven, see
particularly \cite{P,S1,FS,S2,KL,EF,E} and also the survey
\cite{F} for further references. For lattices in other Lie groups
the situation is less clear.  It follows from Gromov's polynomial
growth theorem \cite{Gr1} that any group quasi-isometric to a
nilpotent group is virtually nilpotent, and hence essentially a
lattice in some nilpotent Lie group. However, the quasi-isometry
classification of lattices in nilpotent Lie groups remains an open
problem.

In the case of solvable groups, even less is known.
The main
motivating question is the following:

\begin{conjecture}
\label{conjecture:polycyclic} Let $G$ be a solvable Lie group, and let
$\Gamma$ be a lattice in $G$. Any
finitely generated group $\G'$ quasi-isometric to $\G$ is virtually a
lattice in a (possibly different) solvable Lie group $G'$.
\end{conjecture}

\noindent
{\bf Remarks:}\begin{enumerate}

\item As solvable Lie groups have only cocompact lattices,
studying groups quasi-isometric to lattices in $G$ is equivalent
to studying groups quasi-isometric to $G$.

\item Examples where $G$ and $G'$ need to be different are known.
See \cite{FM3} and Theorem \ref{theorem:abelianbycyclic} below.

\item Conjecture~\ref{conjecture:polycyclic} can be rephrased to
make no reference to connected Lie groups.  In particular, by a
theorem of Mostow, any polycyclic group is virtually a lattice in a
solvable Lie group, and conversely any lattice in a solvable Lie
group is virtually polycyclic \cite{Mo0}.  The conjecture is
equivalent to saying that any finitely generated group
quasi-isometric to a polycyclic group is virtually polycyclic.
This means that being polycyclic is a geometric property.

\item Erschler has shown that a group quasi-isometric to a solvable
  group is not necessarily virtually solvable \cite{D}. Thus, the
  class of virtually solvable groups is not closed under the
  equivalence relation of quasi-isometry. In other words,
  solvability is not a geometric property.

\item Some classes of solvable groups which are not polycyclic are
known to be quasi-isometrically rigid. See particularly the work
of Farb and Mosher on the solvable Baumslag-Solitar groups
\cite{FM1,FM2} as well as later work of Farb-Mosher,
Mosher-Sageev-Whyte and Wortman \cite{FM3,MSW,W}. The methods used
in all of these works depend essentially on topological arguments
based on the explicit structure of singularities of the spaces
studied and cannot apply to polycyclic groups.

\item Shalom has obtained some evidence for the conjecture by
cohomological methods \cite{Sh}.  For example, Shalom shows that
any group quasi-isometric to a polycyclic group has a finite index
subgroup with infinite abelianization. Some of his results have
been further refined by Sauer \cite{Sa}.

\end{enumerate}

Our main results establish Conjecture~\ref{conjecture:polycyclic}
in many cases.  We
believe our techniques provide a method to attack the conjecture.
This is work in progress, some of it
joint with Irine Peng.

\relax From an algebraic point of view, solvable groups are
generally easier to study than semisimple ones, as the algebraic
structure is more easily manipulated. In the present context it is
extremely difficult to see that any algebraic structure is
preserved and so we are forced to work geometrically. For
nilpotent groups the only geometric fact needed is polynomial
volume growth.  For semisimple groups, the key fact for all
approaches is nonpositive curvature. The geometry of solvable
groups is quite difficult to manage, since it involves a mixture
of positive and negative curvature as well as exponential volume
growth.

The simplest non-trivial example for Conjecture
\ref{conjecture:polycyclic} is the $3$-dimensional solvable Lie
group $\SOL$.  This example has received a great deal of
attention. The group $\SOL\cong {\Ra}{\ltimes}\Ra^2$ with $\Ra$
acting on $\Ra^2$ via the diagonal matrix with entries $e^{z/2}$
and $e^{-z/2}$.  As matrices, $\SOL$ can be written as :
\begin{displaymath}
\Sol = \left\{\left. \begin{pmatrix} e^{z/2} & x & 0 \\ 0 & 1 & 0 \\
0 & y & e^{-z/2}
    \end{pmatrix} \right| (x,y,z) \in \Ra^3 \right\}
\end{displaymath}
\noindent The metric $e^{-z}dx^2+e^{z}dy^2 + dz^2$ is a left
invariant metric on $\SOL$.  Any group of the form
${\Za}{\ltimes}_T\Za^2$ for $T \in SL(2,\Za)$ with $|tr(T)|>2$ is
a cocompact lattice in $\SOL$.

The following theorem proves a conjecture of Farb and Mosher and
is one of our main results:

\begin{theorem}
\label{theorem:sol} Let $\G$ be a finitely generated group
quasi-isometric to $\SOL$.  Then $\G$ is virtually a lattice in
$\SOL$.
\end{theorem}

More generally, we can describe the quasi-isometry types of
lattices in many solvable groups.  Here we restrict our attention
to groups of the form $\Ra{\ltimes}\Ra^n$ where the action of
$\Ra$ on $\Ra^n$ is given by powers of an $n$ by $n$ matrix $M$.
The following theorem proves another conjecture of Farb and
Mosher.

\begin{theorem}
\label{theorem:abelianbycyclic} Suppose $M$ is a positive definite
symmetric matrix with no eigenvalues equal to one, and $G=\Ra
\ltimes_M \Ra^n$. If $\Gamma$ is a finitely generated group
quasi-isometric to $G$, then $\Gamma$ is virtually a lattice in
$\Ra \ltimes_{M^{\alpha}} \Ra^n$ for some $\alpha \in \Ra$.
\end{theorem}

\noindent {\bf Remarks:}\begin{enumerate} \item This theorem is
deduced from Theorem \ref{theorem:qiabc} below and a theorem from
the Ph.d. thesis of T.~Dymarz. \item This result is best possible.
All the Lie groups $\Ra \ltimes_{M^{\alpha}} \Ra^n$ for $\alpha \ne 0$
in $\Ra$ are quasi-isometric.
\end{enumerate}

The following is a basic question:
\begin{question}
\label{question:quasilattice} Given a Lie group $G$, is there a
finitely generated group quasi-isometric to $G$?
\end{question}

It is clear that the answer is yes whenever $G$ has a cocompact
lattice. However, many solvable locally compact groups, and in
particular, many solvable Lie groups do not have any lattices.
The simplest examples are groups which are not unimodular.
However, it is possible for Question~\ref{question:quasilattice}
to have an affirmative answer even if $G$ is not unimodular. For
instance, the non-unimodular group solvable group $\left\{
\left.\begin{pmatrix} a & b \\ 0 & a^{-1} \end{pmatrix} \right| a
> 0,
  b \in \Ra \right\}$ acts
simply transitively by isometries on the hyperbolic plane, and
thus is quasi-isometric to the fundamental group of any closed
surface of genus at least $2$.  Thus the answer to
Question~\ref{question:quasilattice} can be subtle. Our methods
give:

\begin{theorem}
\label{theorem:nolattice} Let $G=\Ra{\ltimes}\Ra^2$ be a solvable
Lie group where the $\Ra$ action on $\Ra^2$ is defined by
$z{\cdot}(x,y)=(e^{az}x,e^{-bz}y)$ for $a,b>0$, $a \ne b$.  Then
there is no finitely generated group $\G$ quasi-isometric to $G$.
\end{theorem}

\noindent If $a>0$ and $b<0$, then $G$ admits a left invariant
metric of negative curvature.  The fact that there is no finitely
generated group quasi-isometric to $G$ in this case is a result of
Kleiner \cite{K}, see also \cite{P2}.  Both our methods and
Kleiner's prove similar results for groups of the form
$\Ra{\ltimes}\Ra^n$. Nilpotent Lie groups not quasi-isometric to
any finitely generated group where constructed in \cite{ET}.

In addition our methods yield quasi-isometric rigidity results for
a variety of  solvable groups which are not polycyclic, in
particular the so-called lamplighter groups.   These are the
wreath products $\Za{\wr}F$ where $F$ is a finite group.   The
name lamplighter comes from the description $\Za{\wr}F = F^\Za
\rtimes \Za$ where the $\Za$ action is by a shift. The
subgroup $F^\Za$ is thought of as the states
of a line of lamps, each of which has $|F|$ states.  The
"lamplighter" moves along this line of lamps (the $\Za$ action)
and can change the state of the lamp at her current position. The
Cayley graphs for the generating sets $F \cup \{\pm 1\}$ depend
only on $|F|$, not the structure of $F$.   Furthermore,  $\Za
{\wr} F_1$ and $\Za {\wr} F_2$ are quasi-isometric whenever there
is a $d$ so that $|F_1| = d^s$ and $|F_2|=d^t$ for some $s,t$ in
$\Za$.    The problem of classifying these groups up to
quasi-isometry, and in particular, the question of whether the $2$
and $3$ state lamplighter groups are quasi-isometric, were well
known open problems in the field, see \cite{dlH}.

\begin{theorem}
\label{theorem:lamplighter2} The lamplighter groups $\Za{\wr}F$
and $\Za{\wr}F'$ are quasi-isometric if and only if
there exist
positive integers $d,s,r$ such that $|F|=d^s$ and  $|F'|=d^r$.
\end{theorem}
\noindent For a rigidity theorem for lamplighter groups, see
Theorem~\ref{theorem:lamplighter} below.

To state Theorem \ref{theorem:lamplighter} as well as an analogue
of Theorem \ref{theorem:nolattice} for groups which are not Lie
groups, we need to describe a class of graphs. These are the
Diestel-Leader graphs, $DL(m,n)$, which can be defined as follows:
let $T_1$ and $T_2$ be regular trees of valence $m+1$ and $n+1$.
Choose orientations on the edges of $T_1$ and $T_2$ so each vertex
has $n$ (resp. $m$) edges pointing away from it. This is
equivalent to choosing ends on these trees. We can view these
orientations at defining height functions $f_1$ and $f_2$ on the
trees (the Busemann functions for the chosen ends).   If one
places the point at infinity determining $f_1$ at the top of the
page and the point at infinity determining $f_2$ at the bottom of
the page, then the trees can be drawn as:

\begin{figure}[ht]
\begin{center}
\includegraphics[width=5.5in]{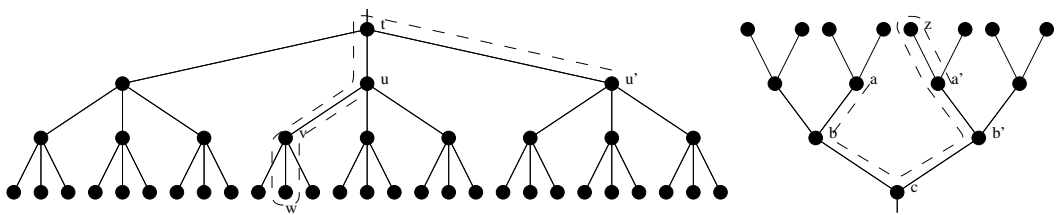}
\caption{The trees for $DL(3,2)$. Figure borrowed from
\cite{PPS}.} \label{fig:1}
\end{center}
\end{figure}

\noindent The graph $DL(m,n)$ is the subset of the product $T_1
\times T_2$ defined by $f_1 + f_2 = 0$. The analogy with the
geometry of $\SOL$ is made clear in section
$\ref{section:solgeom}$.     For $n=m$ the Diestel-Leader graphs
arise as Cayley graphs of lamplighter groups $\Za {\wr} F$ for
$|F|=n$.  This observation was apparently first made by R.Moeller
and P.Neumann \cite{MN} and is described explicitly, from two
slightly different points of view, in \cite{Wo2} and \cite{W}. We
prove the following:

\begin{theorem}
\label{theorem:lamplighter} Let $\G$ be a finitely generated group
quasi-isometric to the lamplighter group $\Za{\wr}F$.  Then there
exists positive integers $d,s,r$ such that $d^s=|F|^r$ and an
isometric, proper, cocompact action of a finite index subgroup of
$\G$ on the Diestel-Leader graph $DL(d,d)$.
\end{theorem}

\noindent{\bf Remark:} The theorem can be reinterpreted as saying
that any group quasi-isometric to $DL(|F|,|F|)$ is virtually a
cocompact lattice in the isometry group of $DL(d,d)$ where $d$ is
as above.
\medskip

In \cite{SW,Wo}, Soardi and Woess ask whether every
homogeneous graph is quasi-isometric to a finitely generated
group. The graph $DL(m,n)$ is easily seen to be homogeneous (i.e. it has
a transitive isometry group).
For $m \neq n$ its isometry group is not unimodular, and hence has
no lattices. Thus there are no obvious groups quasi-isometric to
$DL(m,n)$ in this case.  In fact, we have:

\begin{theorem}
\label{theorem:dl}
There is no finitely generated group
quasi-isometric to the graph $DL(m,n)$ for $m \ne n$.
\end{theorem}
This theorem was conjectured by Diestel and Leader in \cite{DL}, where
the Diestel-Leader graphs were introduced for this purpose. Note that
Theorem~\ref{theorem:dl} can be reinterpreted as the statement that
for $m \ne n$,
there is no finitely generated group quasi-isometric to the isometry
group of $DL(m,n)$.

All our theorems stated above are proved using a new technique,
which we call {\em coarse differentiation}.   Even though
quasi-isometries have no local structure and conventional
derivatives do not make sense, we essentially construct a ``coarse
derivative" that models the large scale behavior of the
quasi-isometry.  This construction is quite different from the
conventional method of passing to the asymptotic cone, see
\S\ref{subsection:coarsedifferentiation} for more discussion.

We now state a theorem that is a well-known consequence of Theorem
\ref{theorem:sol}, Thurston's Geometrization Conjecture and
results in \cite{CC,Gr1,KaL1,KaL2,PW,S1,Ri}.  We state it assuming
that the Geometrization Conjecture is known.

\begin{theorem}
\label{theorem:3manifold} Let $M$ be a compact three manifold
without boundary and $\G$ a finitely generated group.  If $\G$ is
quasi-isometric to the universal cover of $M$, then $\G$ is
virtually the fundamental group of $M'$, also a compact three
manifold without boundary.
\end{theorem}

\section{Quasi-isometries are height respecting} \label{section:qis}

A typical step in the study of quasi-isometric rigidity of groups
is the identification of all quasi-isometries of some space $X$
quasi-isometric to the group, see
\S\ref{section:atinfinity} for a brief explanation.  For us, the
space $X$ is either a solvable Lie group or $DL(m,n)$. In all of
these examples there is a special function $h:X{\rightarrow}\Ra$
which we call the height function and a foliation of $X$ by level
sets of the height function.  We will call a quasi-isometry of any
of these spaces {\em height respecting} if it permutes the height
level sets to within bounded distance (In \cite{FM4}, the term
used is horizontal respecting).

For $\Sol$, the height function is $h(x,y,z)=z$.

\begin{theorem}
\label{theorem:qisol} Any $(K,C)$-quasi-isometry $\varphi$ of
$\SOL$ is within bounded distance of a height respecting  quasi-isometry
$\hat \varphi$.    Furthermore, this distance can be taken
uniform in $(K,C)$ and therefore, in particular,   $\hat
\varphi$ is a $(K',C')$-quasi-isometry where $K',C'$ depend only on
$K$ and $C$.
\end{theorem}

\noindent {\bf Remark:}  In fact, Theorem \ref{theorem:qisol} can
be used to identify the quasi-isometries of $\Sol$ completely.
Possibly after composing with the map
$(x,y,z){\rightarrow}(y,x,-z)$, any height respecting
quasi-isometry (and in particular, any isometry)
is at bounded distance from a quasi-isometry of the form
$(x,y,z){\rightarrow}(f(x),g(y),z)$ where $f$ and $g$ are
bilipschitz functions.  Given a metric space $X$, one defines
$\QI(X)$ to be the group of quasi-isometries of $X$ modulo the
subgroup of those at finite distance from the identity.  The
previous statement can then be taken to mean that
$\QI(\Sol)=\Bilip(\Ra)^2{\ltimes}{\Za/2\Za}$.  This explicit
description was conjectured by Farb and Mosher.

If we take a group of the form $\Ra{\ltimes}\Ra^n$ as in Theorem
\ref{theorem:abelianbycyclic}, we can write coordinates
$(z,\vec{x})$ where $z$ is the coordinate in $\Ra$ and $\vec{x}$
is the coordinate in $\Ra^n$.  Here $h(z,\vec{x})=z$ and level
sets of $h$ are $\Ra^n$ cosets.

\begin{theorem}
\label{theorem:qiabc} Let $X = \Ra{\ltimes}\Ra^n$ be as in Theorem
\ref{theorem:abelianbycyclic}.  Then any $(K,C)$-quasi-isometry
$\varphi$ of $\Ra{\ltimes}\Ra^n$ is within a bounded distance of a
height respecting quasi-isometry $\hat \varphi$. Furthermore,
the bound is uniform in $K$ and $C$.
\end{theorem}

\noindent {\bf Remark:} There is an explicit description of
$QI(\Ra{\ltimes}\Ra^n)$ in this context as well, but it is
somewhat involved so we omit it.

Recall that $DL(m,n)$ is defined as the subset of
$T_{m+1}{\times}T_{n+1}$ where $f_m(x)+f_n(y)=0$ where $f_m$ and
$f_n$ are Busemann functions on $T_m$ and $T_n$ respectively. Here
we simply set $h((x,y))=f_m(x)=-f_n(y)$ which makes sense exactly
on $DL(m,n){\subset}T_{m+1}{\times}T_{n+1}$.  The reader can
verify that the level sets of the height function are orbits for a
subgroup of $\Isom(DL(m,n))$.

\begin{theorem}
\label{theorem:qidl} Any $(K,C)$-quasi-isometry $\varphi$ of
$DL(m,n)$ is within bounded distance from a height respecting
quasi-isometry $\hat \varphi$. Furthermore, the bound is uniform
in $K$ and $C$.
\end{theorem}

\noindent {\bf Remark:} We can reformulate Theorem
\ref{theorem:qidl} in terms similar to those of Theorem
\ref{theorem:qisol}.  Here the group
$\Bilip(\Ra){\times}{\Bilip}(\Ra)$ will be replaced by
$\Bilip(X_m){\times}\Bilip(X_n)$ for $X_m$ (resp. $X_n$) the
complement of a point in the (visual) boundary of $T_{m+1}$ (resp.
$T_{n+1}$).  These can easily be seen to be the $m$-adic and
$n$-adic rationals, respectively.

Note that when $m=n$, this theorem is used to prove Theorem
\ref{theorem:lamplighter} and when $m{\neq}n$ it is used to prove
Theorem \ref{theorem:dl}.  The proofs in these two cases are
somewhat different, the proof in the case $m=n$ being almost
identical to the proof of Theorem \ref{theorem:qisol}.  In the
other case, the argument is complicated by the absence of metric
F\"{o}lner sets, but simplifications also occur since there is no
element in the isometry group that ``flips"
height, see the remarks in \S\ref{subsection:variants}.

There is an analogue of the above results for the case of the
solvable groups which appear in Theorem \ref{theorem:nolattice}.

\section{Geometry of $\SOL$}
\label{section:solgeom}

In this subsection we describe the geometry of $\SOL$ and related
spaces in more detail, with emphasis on the geometric facts used
in our proofs.

The upper half plane model of the hyperbolic
plane $\Ha^2$ is the set $\{(x, \xi) \mid \xi > 0 \}$ with the length element
$ds^2 = \frac{1}{\xi^2} (dx^2 + d\xi^2)$.
If we make the change of
variable $z = \log \xi$, we get $\Ra^2$ with the length element
$ds^2 = dz^2 + e^{-z} dx^2$. This is the {\em log model}
of the hyperbolic plane $\Ha^2$.

The length element of $\Sol$ is:
\begin{displaymath}
ds^2 = dz^2 + e^{-z} dx^2 + e^z dy^2.
\end{displaymath}
Thus planes parallel to the $xz$ plane are hyperbolic planes in
the log model. Planes parallel to the $yz$ plane are {\em
upside-down} hyperbolic planes in the log model. All of these
copies of $\Ha^2$ are isometrically embedded and totally geodesic
.

We will refer to lines parallel to the $x$-axis as $x$-horocycles,
and to lines parallel to the $y$-axis as $y$-horocycles. This
terminology is justified by the fact that each ($x$ or $y$)-horocycle
is indeed a horocycle in the hyperbolic plane which contains it.





We now turn to a discussion of geodesics and quasi-geodesics in
$\SOL$.  Any geodesic in an $\Ha^2$ leaf in $\SOL$ is a geodesic.
There is a special class of geodesics, which we call {\em vertical
geodesics}.  These are the geodesics which are of the form
$\gamma(t)=(x_0,y_0, t)$ or $\gamma(t) = (x_0, y_0, -t)$.
We call the vertical geodesic  {\em upward oriented} in the first case,
and {\em downward oriented} in the second case. In both cases,
this is a unit speed parametrization.
Each vertical geodesic is a geodesic in two hyperbolic planes, the
plane $y=y_0$ and the plane $x=x_0$.

Certain quasi-geodesics in $\Sol$ are easy to describe. Given two
points $(x_0,y_0,t_0)$ and $(x_1,y_1,t_1)$, there is a geodesic
$\gamma_1$ in the hyperbolic plane $y=y_0$ that joins
$(x_0,y_0,t_0)$ to $(x_1,y_0,t_1)$ and a geodesic $\gamma_2$ in
the plane $x=x_1$ that joins $(x_1,y_0,t_1)$ to a $(x_1,y_1,t_1)$.
It is easy to check that the concatenation of $\gamma_1$ and
$\gamma_2$ is a quasi-geodesic. In first matching the $x$
coordinates and then matching the $y$ coordinates, we made a
choice.  It is possible to construct a quasi-geodesic by first
matching the $y$ coordinates and then the $x$ coordinates.  This
immediately shows that any pair of points not contained in a
hyperbolic plane in $\Sol$ can be joined by two distinct
quasi-geodesics which are not close together.  This is an aspect
of positive curvature.  One way to prove that the objects just
constructed are quasi-geodesics is to note the following: The pair
of projections $\pi_1,\pi_2:\Sol{\rightarrow}\Ha^2$ onto the $xt$
and $yt$ coordinate planes can be combined into a quasi-isometric
embedding
$\pi_1{\times}\pi_2:\Sol{\rightarrow}\Ha^2{\times}\Ha^2$.

We state here the simplest version of a key geometric fact used at
various steps in the proof.

\begin{lemma}[Quadrilaterals]
\label{lemma:blocking1} Suppose $p_1,p_2,q_1,q_2 \in \Sol$ and
$\gamma_{ij}: [0,\ell_{ij}] \to \Sol$ are vertical geodesic segments
parametrized by arclength.  Suppose $C > 0$.
Assume that for $i=1,2$, $j=1,2$,
$$d(p_i, \gamma_{ij}(0)) \leq C \qquad \text{ and  } \qquad d(q_i, \gamma_{ij}(\ell_{ij}))
\leq C,$$
so that $\gamma_{ij}$ connects the $C$-neighborhood of $p_i$
to the $C$-neighborhood of $q_j$.  Further assume
that for $i=1,2$ and all $t$,
$d(\gamma_{i1}(t),\gamma_{i2}(t)){\geq} (1/10) t -C$ (so that for
each $i$, the two segments leaving the neighborhood of
$p_i$ diverge right away). Then
there exists $C_1$ depending only on $C$ such that exactly one of the
following holds:
\begin{itemize}
\item[{\rm (a)}] All four $\gamma_{ij}$ are upward oriented,
$p_2$ is within $C_1$ of the
$y$-horocycle passing through $p_1$ and $q_2$ is within
$C_1$ of the $x$-horocycle passing through $\phi(q_1)$.
\item[{\rm (b)}] All four $\gamma_{ij}$ are downward oriented,
$p_2$ is within $C_1$ of the
$x$-horocycle passing through $p_1$ and $q_2$ is within
$C_1$ of the $y$-horocycle passing through $q_1$.
\end{itemize}
\end{lemma}

\noindent We think of $p_1,p_2,q_1$ and $q_2$ as defining a
quadrilateral.  The content of the lemma is that any quadrilateral
has its four "corners" in pairs that lie essentially along
horocycles.
In particular, if we take a quadrilateral with
geodesic segments $\gamma_{ij}$ and with $h(p_1)=h(p_2)$ and
$h(q_1)=h(q_2)$ and map it forward under a $(K,C)$-quasi-isometry $\phi$,
and if we would somehow know that $\phi$ sends
each of the four $\gamma_{ij}$ close to a vertical
geodesic, then Lemma~\ref{lemma:blocking1} would imply that
$\phi$ sends the $p_i$ (resp. $q_i$) to a
pair of points at roughly the same height.

We now define certain useful subsets of $\Sol$. Let $ B(L,\vec{0})
= [-e^L,e^L] \times [-e^L,e^L] \times [-L,L]$. Then
$|B(L,\vec{0})| \approx L e^{2L}$ and $Area(\partial B(L,\vec{0}))
\approx e^{2L}$, so $B(L)$ is a F\"olner set.  We call $B(L,\vec{0})$ a
box of size $L$ centered at the identity.  We define the
box of size $L$ centered at a point $p$ by
$B(L,p)=T_pB(L,\vec{0})$ where $T_p$ is left translation by $p$.
Since left translation is an isometry, $B(L,p)$ is also a
F\"{o}lner set. We frequently omit the center of a box in our
notation and write $B(L)$.

Notice that the top of $B(L)$, meaning the set $[-e^L,e^L] \times
[-e^L,e^L] \times \{L\}$, is not at all square - the sides of this
rectangle are horocyclic segments of  lengths $2e^{2L}$ and $2$
- in other words it is just a small metric neighborhood of a
horocycle.   Similarly, the bottom is also essentially a horocycle
but in the transverse direction.   Further, we can connect the
$1$-neighborhood of  any
point of the top horocycle to the $1$-neighborhood of
any point of the bottom horocycle by
a vertical geodesic segment, and these segments essentially sweep
out the box  $B(L)$.  Thus a box contains an extremely large
number of quadrilaterals.   This picture is even easier to
understand in the Diestel-Leader graphs $DL(n,n)$, where the
boundary of the box is simply the union of the top and bottom
"horocycles", and the vertical geodesics in the box form a
complete bipartite graph between the two.

We remark that for the group in Theorem \ref{theorem:nolattice},
and for the graphs $DL(n,m)$ for $n\neq m$, one has boxes with
essentially the same definition, but these will not be a (metric)
F\"{o}lner set.   A solvable Lie group admits metric F\"{o}lner
sets if and only if it is unimodular.  The same is true of $DL$
graphs.  While the isometry group of a $DL$ graph is always
amenable, the $DL$ graph only has metric F\"{o}lner sets if the
isometry group is unimodular.

\section{On proofs}
\label{section:proofs}

In this section, we give some of the key ideas in the proofs. In
the first two subsections we indicate the key new ideas behind our
proof of Theorem \ref{theorem:qisol}. The first contains
quantative estimates on the behavior of quasi-geodesics. The
second subsection averages this behavior over families of
quasi-geodesics. In \S\ref{subsection:proofsketch} we sketch the
proof of Theorem \ref{theorem:qisol}. Before continuing with
discussion of proofs, we include a discussion of how to axiomatize
the methods of \S\ref{subsec:quasigeodesics} and
\S\ref{subsec:averaging} into a general method of {\em coarse
differentiation} in \S\ref{subsection:coarsedifferentiation}.
Subsection \ref{subsection:variants} briefly discusses the ideas
needed to adapt the proof of Theorem \ref{theorem:qisol} to prove
the other results in Section \ref{section:qis}. In subsection
\S\ref{section:atinfinity}, we discuss deducing results in
\S\ref{section:rigidity} from results in \S\ref{section:qis}.

\subsection{Behavior of quasi-geodesics}
\label{subsec:quasigeodesics}

We begin by discussing some quantative estimates on the behavior
of quasi-geodesic segments in $\SOL$.  Throughout the discussion
we assume $\alpha:[0,r] \to \Sol$ is a $(K,C)$-quasi-geodesic
segment for a fixed choice of $(K,C)$, i.e. $\alpha$ is a
quasi-isometric embedding of $[0,r]$ into $\Sol$.  A
quasi-isometric embedding is a map that satisfies point $(1)$ in
Definition \ref{defn:qi} but not point $(2)$.

\begin{defn}[{\bf $\epsilon$-monotone}]
A quasigeodesic segment $\alpha: [0,r] \to
  \Sol$ is  $\epsilon$-monotone  if for all
  $t_1, t_2 \in [0,r]$ with $h(\alpha(t_1)) =
  h(\alpha(t_2))$ we have $|t_1 - t_2 | < \epsilon r$.
\end{defn}

\begin{figure}[ht]
\begin{center}
\includegraphics[width=1in]{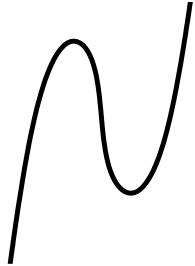}
\caption{A quasigeodesic segment which is not
$\epsilon$-monotone.} \label{fig:not_monotone}
\end{center}
\end{figure}

The following fact about $\varepsilon$-monotone geodesics is an
easy exercise in hyperbolic geometry:

\begin{lemma}[$\epsilon$-monotone is close to vertical]
\label{lemma:monotone} If $\alpha: [0,r] \to \Sol$ is
  $\epsilon$-monotone, then there exists a vertical geodesic segment
  $\lambda$ such that $d(\alpha,\lambda) = O(\epsilon r)$.
\end{lemma}

\noindent{\bf Remark:} The distance $d(\alpha,\lambda)$ is the
Hausdorff distance between the sets and does not depend on
parametrizations.

\begin{lemma}[Subdivision]
\label{lemma:subdivision} Suppose $\alpha: [0,r] \to \Sol$ is a
quasi-geodesic segment which is not $\epsilon$-monotone.  Suppose
$n \GG 1$ (depending on $\epsilon$, $K$, $C$). Then
\begin{displaymath}
\sum_{j=0}^{n-1}\left| h(\alpha(\tfrac{(j+1)r}{n})) -
  h(\alpha(\tfrac{jr}{n})) \right| \ge
\left|h(\alpha(0)) - h(\alpha(r))\right| + \frac{\epsilon r}{ 8
K^2}.
\end{displaymath}
\end{lemma}

\begin{figure}[ht]
\begin{center}
\includegraphics[width=1.0in]{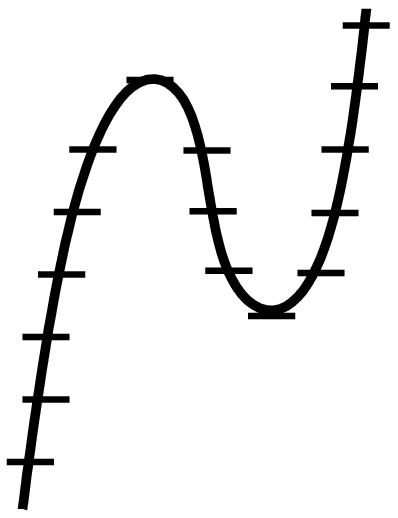}
\caption{Proof of Lemma~\ref{lemma:subdivision}}
\label{fig:subdivision}
\end{center}
\end{figure}

\noindent {\it Outline of Proof.} If $n$ is sufficiently large,
the total variation of the height increases after the subdivision
by a term proportional to $\epsilon$. See
Figure~\ref{fig:subdivision}. \qed\medskip

\noindent {\bf Choosing Scales:} Choose $1 \ll r_0 \ll r_1 \ll
\dots \ll r_M$. In particular, $C \ll r_0$ and $r_{m+1}/r_m > n$.

\begin{lemma}
\label{lemma:scales} Suppose $L \GG r_M$, and suppose $\alpha: [0,
L] \to \Sol$ is a quasi-geodesic segment. For each $m \in [1,M]$,
subdivide $[0,L]$ into $L/r_m$ segments of length $r_m$. Let
$\delta_m(\alpha)$ denote the  fraction  of these segments whose
images are not $\epsilon$-monotone. Then,
\begin{displaymath}
\sum_{m=1}^M \delta_m(\alpha) \le \frac{16K^3}{\epsilon}.
\end{displaymath}
\end{lemma}

\noindent {\it Proof.} By applying Lemma~\ref{lemma:subdivision}
to each non-$\epsilon$-monotone segment on the scale $r_M$, we get
\begin{multline*}
\sum_{j=1}^{L/r_{M-1}} \left|h(\alpha(j r_{M-1})) -
  h(\alpha((j-1)r_{M-1})) \right| \ge \\ \ge
\sum_{j=1}^{L/r_M} \left|
  h(\alpha(j r_M)) -
  h(\alpha((j-1)r_M))\right|+\delta_M(\alpha)
\frac{\epsilon L}{8 K^2}.
\end{multline*}
Doing this again, we get after $M$ iterations,
\begin{multline*}
\sum_{j=1}^{L/r_0} \left|h(\alpha(j r_0)) -
  h(\alpha((j-1)r_0)) \right| \ge \\ \ge
\sum_{j=1}^{L/r_M} \left|
  h(\alpha(j r_M)) -
  h(\alpha((j-1)r_M))\right|+
\frac{\epsilon L}{8 K^2} \sum_{m=1}^M \delta_m(\alpha).
\end{multline*}
But the left-hand-side is bounded from above by the length and so
bounded above by $2K L$. \qed\medskip

\subsection{Averaging}
\label{subsec:averaging}

In this subsection we apply the estimates from above to images of
geodesics under a quasi-isometry of $\SOL$.  The idea is to
average the previous estimates over families of quasi-geodesics.
This results in a coarse analogue of Rademacher's theorem, which
says that a bilipschitz map of $\Ra^n$ is differentiable almost
everywhere, see below for discussion.

{\bf Setup and Notation.}
\begin{itemize}
\item Suppose $\phi: \Sol \to \Sol$ is a $(K,C)$ quasi-isometry.
Without loss of generality, we may assume that $\phi$ is
continuous. \item Let $\gamma: [-L,L] \to \Sol$ be a vertical
geodesic segment parametrized by arclength where $L \gg C$. \item
Let $\overline{\gamma} = \phi \circ \gamma$. Then
$\overline{\gamma}: [-L,L] \to \Sol$ is a quasi-geodesic segment.
\end{itemize}

\noindent It follows from Lemma~\ref{lemma:scales}, that for every
$\theta > 0$ and every geodesic segment $\gamma$, assuming that
$M$ is sufficiently large, there exists $m \in [1,M]$ such that
$\delta_m(\overline{\gamma}) < \theta$. The difficulty is that $m$
may depend on $\gamma$. For $\Sol$, this is overcome as follows:
\medskip

Recall that $ B(L)  = [-e^L,e^L] \times [-e^L,e^L] \times [-L,L]$.
Then $|B(L)| \approx L e^{2L}$ and $Area(\partial B(L)) \approx
e^{2L}$, so $B(L)$ is a F\"olner set. Average the result of
Lemma~\ref{lemma:scales} over $Y_L$, the set of vertical geodesics
in $B(L)$ and let $|Y_L|$ denote the measure/cardinality of $Y_L$.
Changing order, we get:
\begin{figure}
\begin{center}
\includegraphics[width=1.5in]{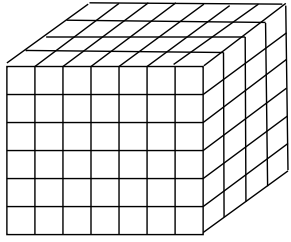}
\caption{The box $B(L)$. } \label{fig:BL}
\end{center}
\end{figure}
\begin{displaymath}
\sum_{m=1}^M \left( \frac{1}{|Y_L|} \sum_{\gamma \in Y_L}
    \delta_m(\overline{\gamma}) \right) \le \frac{32 K^3}{\epsilon}.
\end{displaymath}
Thus, given any $\theta > 0$, (by choosing $M$ sufficiently large)
we can make sure that there exists $1 \le m \le M$ such that
\begin{equation}
\label{eq:avergage:box} \frac{1}{|Y_L|}\sum_{\gamma \in Y_L}
\delta_m(\overline{\gamma}) < \theta.
\end{equation}

{\bf Conclusion.} On the scale $R \equiv {r_m}$, at least
$1-\theta$ fraction of all vertical geodesic segments in $B(L)$
have nearly vertical images under $\phi$. See Figure~\ref{fig:BL}.
\medskip

The difficulty is that, at this point, it may be possible that
some of the (upward oriented) vertical segments in $B(L)$ may have
images which are going up,  and some may have images which are
going down.

We think of the process we have just described as a form of
``coarse differentiation".  For further discussion of this process
and a more general variant on the discussion in the last two
subsections, see subsection
\ref{subsection:coarsedifferentiation}.

\subsection{The scheme of the proof of
  Theorem~\ref{theorem:qisol}.}
\label{subsection:proofsketch}

Roughly speaking, the proof proceeds in the following steps:

\begin{step}
\label{step:I} For all $\theta > 0$ there exists $L_0$ such that
for any box $B(L)$ where $L \geq L_0$, there exists $0 \ll r \ll R
\ll L_0$ such that for the tiling:
\begin{displaymath}
B(L) = \bigsqcup_{i=1}^N B_i(R)
\end{displaymath}
there exists $I \subset \{1,\dots, N\}$ with $|I| \ge (1-\theta)
N$ and for each $i \in I$ there exists a height-respecting map
$\hat{\phi}_i: B_i(R) \to \Sol$ and a subset $U_i \subset B_i(R)$
with $|U_i| \ge (1-\theta)|B_i(R)|$ such that
\begin{displaymath}
d(\phi|_{U_i}, \hat{\phi}_i) = O(r).
\end{displaymath}
\end{step}

Roughly, Step~\ref{step:I} asserts that every sufficiently large
box can be tiled into small boxes, in such a way that for most of
the small boxes $B_i(R)$, the restriction of $\phi$ to $B_i(R)$
agrees, on most of the measure of $B_i(R)$, with a
height-respecting map $\hat{\phi}_i: B_i(R) \to \Sol$. There is no
assertion in Step~\ref{step:I} that the height-respecting maps
$\hat{\phi}_i$ on different small boxes match up to define a
height-respecting map on most of the measure on $B(L)$; the main
difficulty is that some of the $\hat{\phi}_i$ may send the ``up''
direction to the ``down'' direction, while other $\hat{\phi}_i$
may preserve the up direction.

Step~\ref{step:I} follows from a version of
(\ref{eq:avergage:box}) and some geometric arguments using Lemma
\ref{lemma:blocking1}.  The point is that any $\epsilon$-monotone
quasi-geodesic is close to a vertical geodesic by Lemma
\ref{lemma:monotone}.  By the averaging argument in subsection
\ref{subsec:averaging}, we find a scale $R$ at which most segments
have $\epsilon$-monotone image under $\phi$.  More averaging
implies that on most boxes $B_i(R)$ most geodesic segments joining
the top of the box to the bottom of the box have
$\epsilon$-monotone images. We then apply Lemma
\ref{lemma:blocking1} to the images of these geodesics and use
this to show that the map is roughly height preserving on each
$B_i(R)$.  This step also uses the geometric description of
$B_i(R)$ given in the next to last paragraph of
\S\ref{section:solgeom}, i.e. the fact that a box is coarsely a
complete bipartite graphs on nets in the ``top" and ``bottom" of
the box.

\begin{step}
\label{step:II} For all $\theta > 0$ there exists $L_0$ such that
for any box $B(L)$ where $L \geq L_0$, $\exists$ subset $U \subset
B(L)$ with $|U| \ge (1-\theta) |B(L)|$ and a height-respecting map
$\hat{\phi}: B(L) \to \Sol$ such that
$$d(\phi|_U, \hat{\phi}) = O(l),$$
where $l \ll L_0$.
\end{step}

This is the essentially the assertion that the different maps
$\hat{\phi}_i$ from Step~\ref{step:I} are all oriented in the same
way, and can thus be replaced by one standard map $\hat{\phi}:
B(L) \to \Sol$.

Step~\ref{step:II} is the most technical part of the proof. The
problem here derives from exponential volume growth. In Euclidean
space, given a set of almost full measure $U$ in a box, every
point in the box is close to a point in $U$.  This is not true in
$\SOL$ because of exponential volume growth.  Another
manifestation of this difficulty is that $\Sol$ does not have a
Besicovitch covering lemma.  The proof involves using refinements
of Lemma \ref{lemma:blocking1} and further averaging on the image
of $\phi$.

\begin{step}
\label{step:III} The map $\phi$ is $O(L_0)$ from a standard map
$\hat{\phi}$.
\end{step}

This follows from Step~\ref{step:II} and some geometric arguments
using variants of Lemma \ref{lemma:blocking1}. The large constant,
$O(L_0)$, arises because we pass to very large scales to ignore the
sets of small measure that arise in Steps \ref{step:I} and
\ref{step:II}.

\subsection{\bf Remarks on coarse differentiation:}
\label{subsection:coarsedifferentiation}   If a map is
differentiable, then it is locally at sub-linear error from a map
which takes lines to lines. This is roughly the conclusion of the
argument above for the vertical geodesics in $\Sol$, at least on
an appropriately chosen large scale and off of a set of small
measure.  The ideas employed here can be extended to general
metric spaces, by replacing the notion of $\epsilon$-monotone with
a more general notion of $\epsilon$-efficient which we will
describe below. The ideas in our proof are not so different from
the proof(s) of Rademacher's theorem that a bilipschitz map of
$\Ra^n$ is differentiable almost everywhere. In fact, our method
applied to quasi-isometries of $\Ra^n$ gives roughly the same
information as the application of Rademacher's theorem to the
induced bilipschitz map on the asymptotic cone of $\Ra^n$ (which
is again $\Ra^n$). In this context the presence of sets of small
measure can be eliminated by a covering lemma argument.  In the
context of solvable groups, passage to the asymptotic cone is
complicated by the exponential volume growth. The asymptotic cone
for these groups is not locally compact, which makes it difficult
to find useful notions of sets of zero or small measure there.

We now formulate somewhat loosely a more general form of the
``differentiation theorem" given in subsections
\ref{subsec:quasigeodesics} and \ref{subsec:averaging}. Throughout
this subsection $Y$ will be a general metric space, though it may
be most useful to think of $Y$ as a complete, geodesic metric
space. First we generalize the notion of $\epsilon$-monotone.

\begin{defn} \label{definition:eefficient} A quasigeodesic segment
${\alpha}: [0, L] \to Y$ is $\epsilon$-efficient on the scale $r$
if
\begin{displaymath}
\sum_{j=1}^{L/r} d({\alpha}(jr),{\alpha}((j-1)r)) \le (1+\epsilon)
d({\alpha}(L),{\alpha}(0)).
\end{displaymath}
\end{defn}

The fact is that a quasi-geodesic, unless it is a $(1+\epsilon)$
quasi-geodesic, fails to be $\epsilon$-efficient at some scale some
fraction of the time.  The observation embedded in subsection
\ref{subsec:quasigeodesics} is that this cannot happen everywhere
on all scales and in fact cannot happen too often on too many
scales.

\begin{figure}[ht]
\begin{center}
\includegraphics[width=3in]{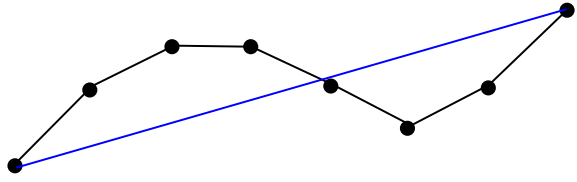}
\caption{The definition of $\epsilon$-efficient.} \label{fig:6}
\end{center}
\end{figure}

With this definition, the following variant on Lemma
\ref{lemma:subdivision} becomes a tautology.

\begin{lemma}[Subdivision II]
\label{lemma:subdivision2} Given $\epsilon>0$, there exist $r \gg
C$ and $n\gg 1$ (depending on $K,C$ and $\epsilon$) such that any
$(K,C)$-quasi-geodesic segment $\alpha: [0,r] \to X$ which is not
$\epsilon$-efficient on scale $\frac{r}{n}$ we have:
\begin{displaymath}
\sum_{j=0}^{n-1}
d(\alpha(\tfrac{(j+1)r}{n}),\alpha(\tfrac{jr}{n})) \ge
d(\alpha(0),\alpha(r)) + \frac{\epsilon r}{2K}.
\end{displaymath}
\end{lemma}

We now state a variant of Lemma \ref{lemma:scales} whose proof is
verbatim the proof of that lemma.

\noindent {\bf Choosing Scales:} Choose $1 \ll r_0 \ll r_1 \ll
\dots \ll r_M$. In particular, $C \ll r_0$ and $r_{m+1}/r_m > n$.

\begin{lemma}
\label{lemma:scales2} Suppose $L \GG r_M$, and suppose $\alpha:
[0, L] \to X$ is a quasi-geodesic segment. For each $m \in [1,M]$,
subdivide $[0,L]$ into $L/r_m$ segments of length $r_m$. Let
$\delta_m(\alpha)$ denote the  fraction  of these segments whose
images are not $\epsilon$-efficient on scale ${r_{m-1}}$. Then,
\begin{displaymath}
\sum_{m=1}^M \delta_m(\alpha) \le \frac{4K^2}{\epsilon}.
\end{displaymath}
\end{lemma}

Let $X$ be a geodesic metric space.  Coarse differentiation
amounts to the following easy lemma.

\begin{lemma}[Coarse Differentiation]
\label{lemma:coarsedifferentiation} Let $\phi:X{\rightarrow}Y$ be
a $(K,C)$-quasi-isometry.  For all $\theta > 0$ there exists $L_0
\gg 1$ such that for any $L>L_0$ and any family $\mathcal F$ of
geodesics of length $L$ in $X$, there exist scales $r,R$ with $C
\ll r \ll R \ll L_0$ such that if we divide each geodesic in
$\mathcal{F}$ into subsegments of length $R$, then at least
$(1-\theta)$ fraction of these subsegments have images which are
$\epsilon$-efficient at scale $r$.
\end{lemma}

This lemma and its variants seem likely to be useful in other
settings. In fact, the lemma holds only assuming that $\phi$ is
coarsely lipschitz. A map $\phi:X{\rightarrow}Y$ is a {\em $(K,C)$
coarsely lipschitz} if
$d_Y(\phi(x_1),\phi(x_2)){\leq}Kd_X(x_1,x_2)+C$. We now describe
the relation to taking derivatives and also to the process of
taking a ``derivative at infinity" of a quasi-isometry by passing
to asymptotic cones.

We first discuss the case of maps $\Ra^n{\rightarrow}\Ra^n$.
Suppose $\phi:\Ra^n{\rightarrow}\Ra^n$ is a quasi-isometry.
Suppose one chooses a net $N$ on the unit
circle and takes $\mathcal{F}$ to be the set of all lines of
length $L$ in a large box, whose direction vector is in $N$.
Lemma~\ref{lemma:coarsedifferentiation} applied to $\mathcal{F}$
then states that most
of these lines, on the appropriate scale, map under $\phi$ close to straight
lines, which implies that the map $\phi$ (in a suitable box) can be
approximated by an affine map. Thus, in this context,
Lemma~\ref{lemma:coarsedifferentiation} is indeed analogous to
differentiation (or producing points of differentiability).

An alternative approach for analyzing quasi-isometries
$\phi:\Ra^n{\rightarrow}\Ra^n$ is to pass to the asymptotic cone to
obtain a bilipschitz map $\tilde \phi:\Ra^n{\rightarrow}\Ra^n$ and
then apply Rademacher's theorem to $\tilde{\phi}$. If one attempts to
pull the information this yields back to $\phi$ one gets statements
that are similar to those one would obtain directly using Lemma
\ref{lemma:coarsedifferentiation}. This is not surprising, since
averaging arguments like those used in the proof of Lemma
\ref{lemma:coarsedifferentiation} are implicit in the proofs of
Rademacher's theorem.

Passing to the asymptotic cone has obvious advantages because it
allows one to replace a $(K,C)$ quasi-isometry from $X$ to $Y$ with
a $(K,0)$-quasi-isometry (i.e. a bilipschitz map) from the asymptotic
cone of $X$ to the asymptotic cone of $Y$. One can then try to use
analytic techniques to study the
bilipshitz maps. However, a major difficulty which occurs is that the
asymptotic cones are typically not locally compact and notions of
measure and averaging on such spaces are not clear.  This difficulty
arises as soon as one has exponential volume growth. In particular
it is not clear if there is a useful version of Rademacher's theorem
for the asymptotic cones of the spaces which we consider in this
paper.

The main advantage of Lemma~\ref{lemma:coarsedifferentiation} compared
to the asymptotic cone approach is that the averaging is done on the
(typically locally compact)
space $X$, i.e. the domain of the quasi-isometry $\phi$. In other
words, we construct a ``coarse derivative'' without first
passing to a limit to get rid of the additive constant. In particular,
the information we obtain about $\Sol$ and other solvable groups by
coarse differentiation is not easily extracted by passage to the
asymptotic cone.

We remark again that Lemma~\ref{lemma:coarsedifferentiation}
applies to any quasi-isometric embedding (or any uniform
embedding) between any two metric spaces $X$ and $Y$. However its
usefulness clearly depends on the situation.

\subsection{Remarks on Theorems \ref{theorem:qiabc} and
\ref{theorem:qidl}} \label{subsection:variants}

The proof of Theorem \ref{theorem:qiabc} is quite similar to the
proof of Theorem \ref{theorem:qisol} but becomes much more
involved technically in a few places, particularly at Step
\ref{step:II}.

The use of F\"{o}lner sets in the proof of Theorem
\ref{theorem:qisol} might make it surprising that similar
techniques apply to prove  Theorems \ref{theorem:nolattice},
\ref{theorem:dl}
 and \ref{theorem:qidl}. As remarked in
Section \ref{section:solgeom}, it is well known that there are no
(metric) F\"{o}lner sets for $DL(m,n)$ when $m{\neq}n$ or for
non-unimodular solvable Lie groups. The key is to use a notion of
weighted amenability and weighted averaging that is similar to the
one used in \cite{BLPS}. In our setting this arises quite
naturally. We are averaging over the set of geodesics in a box.
The asymmetry of boxes in this context implies that points near
the ``top" of the box are on more geodesics than points near the
``bottom".  This reweighting process introduces a new measure
which is not, a priori, quasi-invariant under quasi-isometries. It
is easy to see that the standard volume is quasi-invariant under
quasi-isometries.  The new measure is a reweighting of the
standard volume by a factor depending only on the height. Using
Lemma \ref{lemma:blocking1} and its variants to see that height is
coarsely preserved allows us to also conclude that this new
measure is coarsely preserved. The argument at Step \ref{step:II}
then simplifies dramatically, since we can show that no
quasi-isometry can ``flip" the orientation of a box.

\subsection{Deduction of rigidity results}
\label{section:atinfinity}

In our setting, the deduction of rigidity results from the
classification of quasi-isometries follows a fairly standard
outline that is similar to one used for semisimple groups as well
as for certain solvable groups in \cite{FM2,FM3,MSW}.  As this is
standard, we will say relatively little about it. Some of these
ideas go back to Mostow's original proof of Mostow rigidity
\cite{Mo1,Mo2} and have been developed further by many authors.

Given a group $\G$ any element of $\g$ in $\G$ acts on $\G$ by
isometries by left multiplication $L_{\g}$.  If $X$ is a metric
space and $\phi:\G{\rightarrow}X$ is a quasi-isometry, we can
conjugate each $L_{\g}$ to a  self quasi-isometry
$\phi{\circ}L_{\g}{\circ}\phi{\inv}$ of $X$. This induces a
homomorphism of $\Phi:\G\to \QI(X)$. Here $\QI(X)$ is the group of
quasi-isometries of $X$ modulo the subgroup of quasi-isometries a
bounded distance from the identity.  The approach we follow is to
use $\Phi$ to define an action of $\G$ on a ``boundary at
infinity" of the space $X$.  All theorems are then proven by
studying the dynamics of this ``action at infinity."  We are
ignoring many important technical points here, such as why $\Phi$
has finite kernel and why $\QI(X)$ acts on either $X$ or the
boundary at infinity of $X$.

The deduction of Theorem \ref{theorem:sol} from Theorem
\ref{theorem:qisol} was known to Farb and Mosher \cite{FM2,FM4}.
The action at infinity is studied using a variant of a theorem of
Hinkannen due to Farb and Mosher \cite{H,FM2,FM4}. In the context
of Theorem \ref{theorem:abelianbycyclic}, we deduce the result
from Theorem \ref{theorem:qiabc} using results from the
dissertation of Tullia Dymarz \cite{Dy}.  These are variants and
extensions of the results of Tukia in \cite{Tu}.

In the context of Diestel-Leader graphs the argument is somewhat
different than in the previous cases.  In this context we use
\cite[Theorem $7$]{MSW} to understand the dynamics at infinity.
While this result was motivated by analogy with the results
discussed above, its proof is quite different, and uses topology
in place of analysis.  The use of \cite[Theorem $7$]{MSW} is
precisely the step in the proof of Theorem
\ref{theorem:lamplighter} where we might need to replace
$DL(|F|,|F|)$ with $DL(d,d)$ where $d$ and $F$ are powers of a
common integer. A similar argument occurs in the proof of Theorem
\ref{theorem:dl}.

\noindent  Department of Mathematics, University of Chicago, 5734
S. University Avenue, Chicago, Illinois 60637.\

\noindent Department of Mathematics, Indiana University, Rawles
Hall, Bloomington, IN, 47405.\

\noindent Department of Mathematics, Statistics, \& Computer
Science, 322 Science \& Engineering Offices (M/C 249), 851 S.
Morgan Street Chicago, IL 60607-7045.\

\end{document}